% stidelft.tex
%
%  This is plain TeX
%
% A.Magnus, summer 1994 , revised, Feb. 1995
%
%%
%% Here is a character listing to check to be sure that no
%% unwanted translations took place within the bowels of the net.
%%
%%   Upper-case     A B C D E F G H I J K L M N O P Q R S T U V W X Y Z
%%   Lower-case     a b c d e f g h i j k l m n o p q r s t u v w x y z
%%   Digits         0 1 2 3 4 5 6 7 8 9
%%   Exclamation    !     Double quote   "     Hash (number)  #
%%   Dollar         $     Percent        %     Ampersand      &
%%   Acute accent   '     Left paren     (     Right paren    )
%%   Asterisk       *     Plus           +     Comma          ,
%%   Minus          -     Point          .     Solidus,slash  /
%%   Colon          :     Semicolon      ;     Less than      <
%%   Equals         =     Greater than   >     Question mark  ?
%%   Commercial at  @     Left bracket   [     Backslash      \
%%   Right bracket  ]   Circumflex,caret ^     Underscore     _
%%   Grave accent   `     Left brace     {     Vertical bar   |
%%   Right brace    }     Tilde          ~
%%

 \magnification=1250
%\vsize=10.0 truein
% extrait de manmac.tex
% Macros for The TeXbook
 
\catcode`@=11 % borrow the private macros of PLAIN (with care)
 
\font\eightrm=cmr10 scaled 800
 
\font\sixrm=cmr7 scaled 800
\font\sixbf=cmbx7 scaled 800
 
\font\eighti=cmmi10 scaled 800
\font\sixi=cmmi7 scaled 800
 
\font\eightsy=cmsy10 scaled 800
\font\sixsy=cmsy7 scaled 800

\font\eightbf=cmbx10 scaled 800
 
\font\eighttt=cmtt10 scaled 800
 
\font\eightsl=cmsl10 scaled 800
 
\font\eightit=cmti10 scaled 800
 
\font\fiverm=cmr5 scaled 800
\font\fivei=cmmi5 scaled 800
\font\fivesy=cmsy5 scaled 800
 
\font\tenux=cmex10 scaled 800
 
\def\eightpoint{\def\rm{\fam0\eightrm}%
  \textfont0=\eightrm \scriptfont0=\sixrm \scriptscriptfont0=\fiverm
  \textfont1=\eighti \scriptfont1=\sixi \scriptscriptfont1=\fivei
  \textfont2=\eightsy \scriptfont2=\sixsy \scriptscriptfont2=\fivesy
  \textfont3=\tenux \scriptfont3=\tenux \scriptscriptfont3=\tenux
  \def\it{\fam\itfam\eightit}%
  \textfont\itfam=\eightit
  \def\sl{\fam\slfam\eightsl}%
  \textfont\slfam=\eightsl
  \def\bf{\fam\bffam\eightbf}%
  \textfont\bffam=\eightbf \scriptfont\bffam=\sixbf
   \scriptscriptfont\bffam=\fivebf
  \def\tt{\fam\ttfam\eighttt}%
  \textfont\ttfam=\eighttt
  \normalbaselineskip=9pt
  \def\MF{{\manual opqr}\-{\manual stuq}}%
  \let\sc=\sixrm
  \let\big=\eightbig
  \setbox\strutbox=\hbox{\vrule height7pt depth2pt width\z@}%
  \normalbaselines\rm}
 
\catcode`@=12

\def\D{{\cal D}}

\def\M{{\cal M}}
 
\def\strutt{\vrule width 0 pt height 9.5 pt depth 4.0pt}
\def\dfrac#1#2{{\displaystyle{#1\over#2}}}
                        \tolerance=999
\def\QED{$\sqcup\kern-.65em\sqcap$}
\overfullrule=0pt

\input epsf.tex    \epsfverbosetrue
 
%\nopagenumbers
\vglue 25pt
 
\centerline{\bf Special non uniform lattice ($snul$) orthogonal polynomials}
\centerline{\bf on discrete dense sets of points.}
 
\bigskip
\centerline{Alphonse P. Magnus}
 
\centerline{ Institut Math\'ematique,  Universit\'e Catholique de Louvain}
 
\centerline{ Chemin du Cyclotron 2}
 
\centerline { B-1348 Louvain-la-Neuve}
 
\centerline{ Belgium}
 
\centerline{ e-mail:  {\tt magnus@anma.ucl.ac.be}}
 
\bigskip\bigskip
 
\noindent{\sl Keywords:\/} Orthogonal polynomials, difference
    operators.
 
\bigskip\bigskip
\rightline{\sl ``Il n'est pas n\'ecessaire d'esp\'erer pour entreprendre,}
\rightline{\sl ni de r\'eussir pour pers\'ev\'erer.''\qquad}
 
\rightline{(\sl {Begin, even without hope; Proceed, even without success.})}
 
\rightline{\sl William of Orange (``William the Silent''), murdered in
   Delft in 1584.}
 
\bigskip
 
{\bf Abstract.\/} Difference calculus compatible with polynomials
(i.e., such that the divided difference operator of first order applied
to any polynomial must yield a polynomial of lower degree) can only be
made on special lattices well known in contemporary $q-$calculus.
Orthogonal polynomials satisfying difference relations on such lattices
are presented. In particular, lattices which are dense on intervals
($|q|=1$) are considered.
 
\bigskip
 
{\bf 1. Introduction.}
 
Many works have been devoted to orthogonal polynomials satisfying
remarkable differential or difference relations.
 
For instance, the classical orthogonal polynomials are characterized
by the existence of a differential relation of the form
$$ W(x)p'_n(x) = \omega_n(x) p_n(x) +\vartheta_n p_{n-1}(x), \eqno(1)$$
where $W$ is a fixed polynomial of degree $\le 2$,
      $\omega_n$ is a polynomial of degree $\le 1$,and where
$\vartheta_n$ is constant (in $x$) (cf. [1], p.8).
 
This differential relation interacts  most efficiently with the
recurrence relation
$$ a_{n+1}p_{n+1}(x) = (x-b_n)p_n(x) +a_n p_{n-1}(x) \eqno(2)$$
in the production of various useful identities.
 
If we accept higher degrees polynomials in (1):
$$ W(x)p'_n(x) = \Omega_n(x)p_n(x) -a_n\Theta_n(x)p_{n-1}(x),\eqno(3)$$
with $W,\Omega_n$ and $\Theta_n$ polynomials of degrees $\le s+2,s+1$ and
$s\ge 0$, we get the {\it semi-classical\/} class, studied by
Laguerre (the notations of (3) are [almost] Laguerre's ones [19]),
Hendriksen and van Rossum [10], and Maroni [22] who coined this
name.
See also  [4,5] for determination of the relevant measure.
 
If we try to extend (3) to a difference operator of first order, we
expect to see the derivative $p'_n(x)$ replaced by some combination
of $p_n(y(s))$ and $p_n(y(s+1))$, where $y(s)$ and $y(s+1)$ are two
consecutive points on a  lattice associated to the difference operator.
 
We will explore here the extension of the semi-classical property (3) to
the remarkable difference operators and the corresponding nonuniform
lattices studied by many people recently [2,6,12,16,17,21,
24,25,29,30].
In particular, those lattices which happen to fill densely an
interval will be considered with special care.
\bigskip
 
{\bf 2. The difference operator and the related lattices.}

\bigskip
 
We consider here a first-order difference operator involving the
values of a function at two points. For each $x$, let $\varphi_1(x)$
and $\varphi_2(x)$ be these still unknown points. The first-order
{\it divided difference\/}
operator at $x$ is
 $$(\D f)(x) =\dfrac{f(\varphi_2(x))-f(\varphi_1(x))}
                    {\varphi_2(x)-\varphi_1(x)}\;. \eqno(4) $$
If we impose the condition
 that $\D f$ is a polynomial of degree $n-1$ if $f$ has degree $n$,
then $\varphi_1(x)$ and $\varphi_2(x)$ must be the two roots in $y$ of
a quadratic equation
 $$Ay^2+2Bxy+Cx^2+2Dy+2Ex+F=0. \eqno(5) $$
  (see [6,12, 21]).
 
Indeed, applying $\D$ to $f(x)=x^2$ and $f(x)=x^3$ readily yields that
$\varphi_1+\varphi_2$ and
$\varphi_1^2+\varphi_1\varphi_2+\varphi_2^2$ must be polynomials of
degrees 1 and 2, which implies (5). Conversely, if (5) holds, any
symmetric polynomial in $\varphi_1$ and $\varphi_2$ is a polynomial
in $x$.
 
Let us figure the conic (5) and one of its parametric representations
$\{x(s),y(s)\}$ such that $y(s)$ and $y(s+1)$ appear naturally as the
two ordinates associated to the abscissa $x=x(s)$:
one starts from some point $\{x_1=x(s_1),y_1=y(s_1)\}$ on the conic,
and one looks for the points $\{x_k=x(s_1+k-1),y_k=y(s_1+k-1)\}$,
$k=1,2,\ldots$
 
\epsfxsize=477bp
\epsfbox[0 150 500 400]{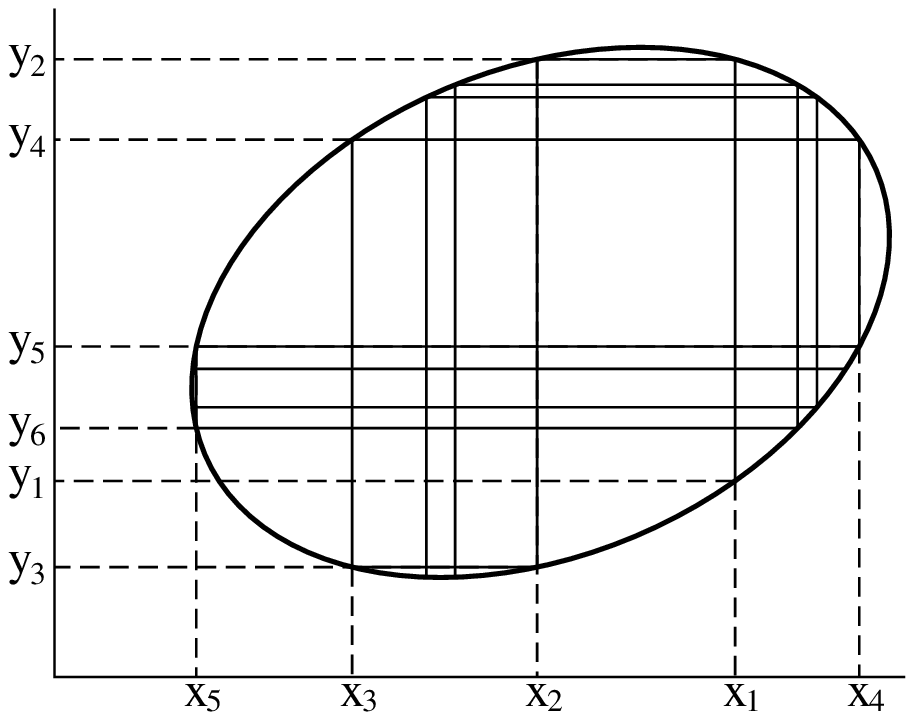}
\medskip
\centerline{Figure 1. $x-$ and $y-$ lattices.}
 
\bigskip
\bigskip
 
To achieve this, let us consider the familiar parametric representations
in the two following cases:
 
\item{1.} The conic (5) has a center $B^2-AC\ne 0$. With the center
coordinates
  $x_c=(AE-BD)/(B^2-AC)$ and
  $y_c=(CD-BE)/(B^2-AC)$, one has
$A(y-y_c)^2+2B(x-x_c)(y-y_c)+C(x-x_c)^2+\widetilde F=0$, with
$\widetilde F=F-Ay_c^2-2Bx_cy_c-Cx_c^2=F+Dy_c+Ex_c =
F+(CD^2-2BDE+AE^2)/(B^2-AC)$,
$$ x=x(s)=x_c+\zeta \sqrt{A}(q^s+q^{-s})\;,
   y=y(s)=y_c+\zeta \sqrt{C}(q^{s-1/2}+q^{-s+1/2})\;,\eqno(6) $$
is a valid parametric representation of (5), (if $AC\ne 0$), where
\hfill\break
$\zeta^2 = \widetilde F/(4(B^2-AC))$ and
$$ q^{1/2}+q^{-1/2}=-\dfrac{2B}{\sqrt{AC}}\; {\rm \ i.e.,\ }
  q+q^{-1}=\dfrac{4B^2}{AC}-2. \eqno(7)$$
Indeed, $x(s)=x(-s)$ is kept by the transformation
$s\leftrightarrow -s$ in (6) but $y$ becomes $y(s+1)$.
\hfill\break
If $AC=0$, the scheme of Fig.~1  and Fig.~2
does not work: horizontal and/or
vertical lines do not meet the conic in two points any more.
 
\smallskip
 
\epsfxsize=477bp
\epsfbox[50  150 527 475]{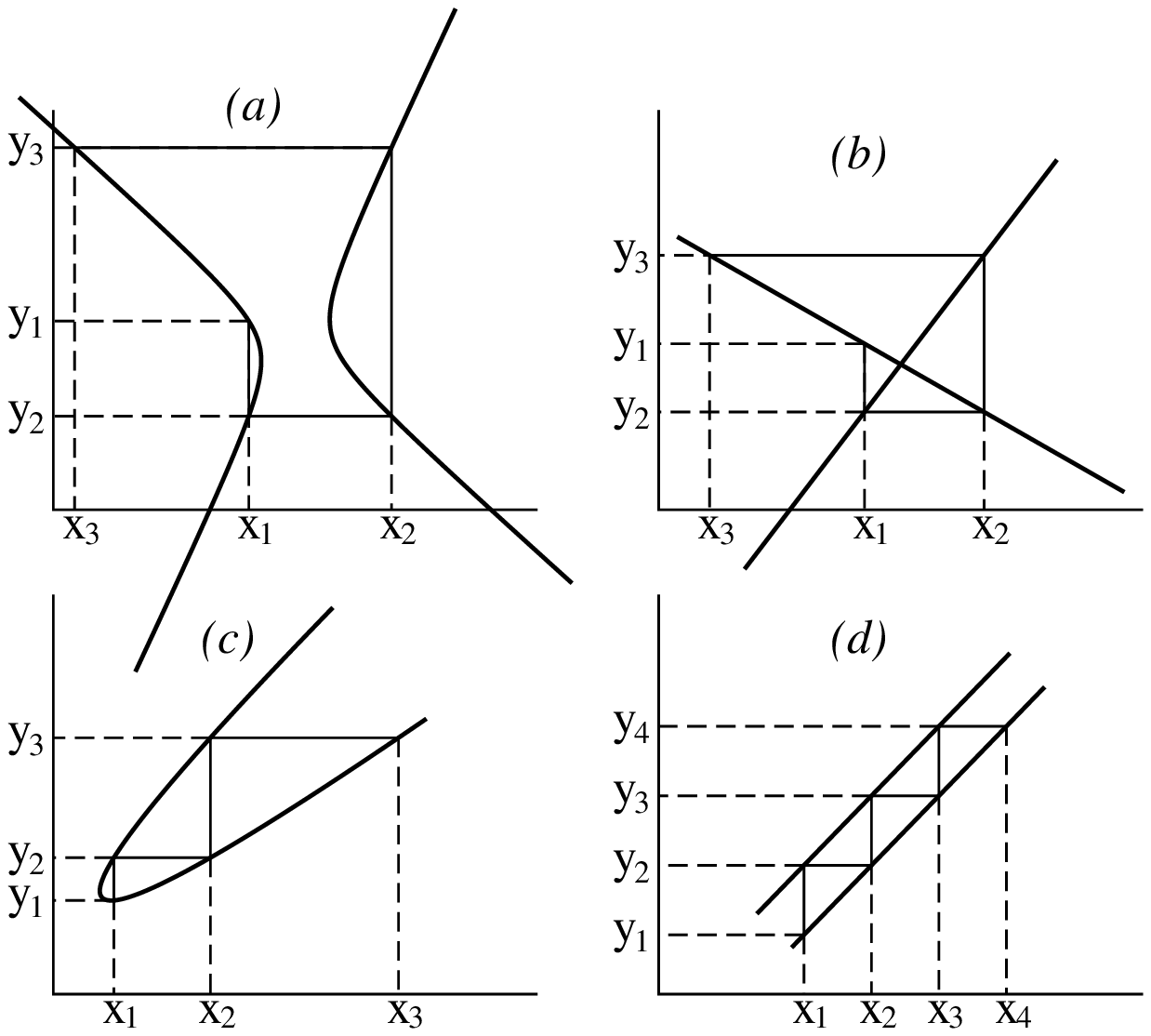}
\medskip
\centerline{Figure 2. Other kinds of lattices.}
 
%\bigskip
\bigskip
\noindent
The generic (also rightly called {\it hyperbolic\/}) case $|q|\ne 1$
gives a hyperbola (Fig.~2(a)).
  Remark that the asymptotes are given by
($x$ and $y$ large) $y\sim (C/A)^{1/2} q^{\pm 1/2} x$:
\hfill\break
{\centerline{$q$ {\it is the ratio of the slopes of the asymptotes of the conic.}}}
\hfill\break
If $\widetilde F=0$, one finds $x-x_c=X\sqrt{A} q^s$,
                               $y-y_c=X\sqrt{C} q^{s\pm 1/2}$,
the ``old'' $q-$lattice (Fig.~2(b)).
 
\item{2.} The conic (5) has no center, $B^2-AC =  0$. Then,
$$\eqalign{
   x=x(s)&= \sqrt{A}\left\{
   \dfrac{D^2-AF}{2A(D\sqrt{C}+E\sqrt{A})}
 -2  \dfrac{D\sqrt{C}+E\sqrt{A}}{AC} s^2 \right\}\;,  \cr
   y=y(s)&= \sqrt{C}\left\{
   \dfrac{E^2-CF}{2C(D\sqrt{C}+E\sqrt{A})}
 -2  \dfrac{D\sqrt{C}+E\sqrt{A}}{AC} (s-1/2)^2 \right\}\;, \cr}
$$
found directly to satisfy (5) in the parabolic case $B^2=AC$, or by
taking the limit of (6) when $q\to 1$: let $q=\exp(\varepsilon)$, then,
$q^s=1+s\varepsilon+s^2\varepsilon^2/2+\cdots$,
$B=-\sqrt{AC}(q^{1/2}+q^{-1/2})/2 =
   -\sqrt{AC}(1+\varepsilon^2/8+\cdots)$, etc. (Fig.~2(c)).
\hfill\break
If
$D\sqrt{C}+E\sqrt{A} =0$, one redefines $s$ through the translation
$$s_{old} = \lim_{D\sqrt{C}+E\sqrt{A}\to 0}\pm\left\{ s_{new} +
   \dfrac{\sqrt{(D^2-AF)C}}{2(D\sqrt{C}+E\sqrt{A})}\right\}\;,$$
to get
$$x=-2\sqrt{\dfrac{D^2-AF}{AC}} s\;,
  y=\dfrac{E}{\sqrt{AC}}-2\dfrac{\sqrt{D^2-AF}}{A} (s\pm 1/2)\;, $$
the simplest lattice (this one is uniform (Fig.~2(d))).
 
\smallskip
 
\noindent
These points form one of the special non uniform lattices ($snul$)
 I-VI of
 [24,25].
If the quadratic equation (5) describes an ellipse, as in the figure, this
lattice fills densely an interval, unless it is finite (periodic:
$q^N=1$).
\bigskip
 
{\bf 3. Semi-classical orthogonal polynomials on $snul$s.}

\bigskip

Semi-classical $snul$ orthogonal polynomials may be defined through a
$\D-$difference equation of the form
$$W(x) (\D S)(x)=2V(x)(\M S)(x)+U(x)  \eqno(8)$$
 for the
Stieltjes function
$$ S(x)=\int_{{\rm Supp.\ }\mu}(x-t)^{-1}d\mu(t)=
\sum_0^\infty \mu_k/x^{k+1}$$
(Stieltjes transform of the orthogonality measure $d\mu$)
where $W$, $V$ and $U$ are polynomials and $\M$ is the arithmetic
mean operator:
$$(\M f)(x)=(f(\varphi_1(x))+f(\varphi_2(x)))/2.$$
 
 If $\mu$ has a jump (pole of $S$) at some $y_n$,
it must have jumps at the other points $y_{n\pm 1}$, $y_{n\pm 2}$,\dots
of the corresponding lattice.
 
It is possible to recover a difference relation extending (3) to
the present difference operators, as attempted and
{\it min of meer\/} (less or more) achieved in the paper [21],
see here an attempt to achieve the converse:
 
\medskip
 
\noindent{\bf Theorem.\/} {\sl
 Orthogonal polynomials are $\D-$semi-classical,
i.e., their Stieltjes function satisfies an equation $(8)$ with
polynomials $W$, $V$, and $U$, if each $p_n$ satisfies a linear
first-order difference relation connecting $p_n$ with $p_{n-1}$
 $$ W_n(x) (\D p_n)(x) = \Omega_n(x) (\M p_n)(x) -a_n\Theta_n(x)
  (\M p_{n-1})(x), \eqno(9) $$
where $W_n$, $\Omega_n$, and $\Theta_n$ are polynomials of fixed
$($independent of $n)$ degrees.}
 
\medskip
 
Indeed, let $y_1$ and $y_2$ be the ordinates corresponding to $x=x_1$,
then (9) is a linear relation involving $p_n(y_1)$, $p_n(y_2)$,
$p_{n-1}(y_1)$, and $p_{n-1}(y_2)$:
 $$ A_n p_n(y_1)+B_n p_n(y_2) +C_n p_{n-1}(y_1) +D_n p_{n-1}(y_2)=0,
       \eqno(10) $$
where $A_n$, $B_n$, $C_n$, and $D_n$ are rational functions of fixed
degrees of $x=x_1$  and $y_1$ ($y_2$ may be replaced by
$-2(Bx+D)/A-y_1$ in $A_n=-W_n(x)/(y_2-y_1)-\Omega_n(x)/2$ etc.).
  We come to $n+1$:
 $$ A_{n+1} p_{n+1}(y_1)+B_{n+1} p_{n+1}(y_2) +C_{n+1} p_n(y_1)
  +D_{n+1} p_n(y_2)=0, $$
and use the recurrence relation (2):
$$\displaylines{
  \left[(y_1-b_n)A_{n+1}+a_{n+1}C_{n+1}\right]p_n(y_1) +
  \left[(y_2-b_n)B_{n+1}+a_{n+1}D_{n+1}\right]p_n(y_2) - \cr
    \hfill  a_n A_{n+1}p_{n-1}(y_1) -
  a_n B_{n+1}p_{n-1}(y_2) =0,  \qquad(11) \cr}$$
which is also a linear relation involving $p_n(y_1)$, $p_n(y_2)$,
$p_{n-1}(y_1)$, and $p_{n-1}(y_2)$!
\hfill\break
 
\item{1.} Either (10) and (11) are dependent for each $n$, then
$A_n=K_n[(y_1-b_n)A_{n+1}+a_{n+1}C_{n+1}]$ and $C_n=K_n a_nA_{n+1}$,
i.e., the recurrence relation
$A_n=K_n[(y_1-b_n)A_{n+1}-a_{n+1}^2K_{n+1}A_{n+2}]$
which looks somewhat like (2): actually,\hfill\break
 $K_1\ldots K_{n-1}
a_1\ldots a_{n-1}A_n$ would be a solution of (2), and this is
incompatible with the requirement that $A_n$ keeps a finite degree.
So,
 
\item{2.} or (10) and (11) are independent for some $n$. Then, it is
possible to extract $p_n(y_1)$ and $p_n(y_2)$ in terms of
                    $p_{n-1}(y_1)$ and $p_{n-1}(y_2)$
:
$$ p_n(y_1)=X_np_{n-1}(y_1)+Y_np_{n-1}(y_2)\;,\quad
   p_n(y_2)=Z_np_{n-1}(y_1)+U_np_{n-1}(y_2),$$
where $X_n$, $Y_n$, $Z_n$, and $U_n$ are again rational functions
of fixed degrees. We get the relation for $n+1$ using (2):
$$ \eqalign{
   \left[\matrix{p_{n+1}(y_1)\cr p_{n+1}(y_2)\cr }\right] &=
   \dfrac1{a_{n+1}}
   \left[\matrix{ y_1-b_n & 0 \cr 0 & y_2-b_n \cr}\right]
   \left[\matrix{p_n(y_1)\cr p_n(y_2)\cr }\right] -
   \dfrac{a_n}{a_{n+1}}
   \left[\matrix{p_{n-1}(y_1)\cr p_{n-1}(y_2)\cr }\right]  \cr
   &= \left\{
   \dfrac1{a_{n+1}}
   \left[\matrix{ y_1-b_n & 0 \cr 0 & y_2-b_n \cr}\right]
 - \dfrac{a_n}{a_{n+1}}
   \left[\matrix{X_n & Y_n \cr Z_n & U_n \cr}\right]^{-1}
   \right\}
   \left[\matrix{p_n(y_1)\cr p_n(y_2)\cr }\right] \;, \cr
       }
$$
showing that we can construct similar relations for $n+1$, $n+2$, etc.
as far as $x$ is not a zero of the determinants
$\delta_n=X_nU_n-Z_nY_n$, etc.:
$$
   \left[\matrix{p_{n+1}(y_1)\cr p_{n+1}(y_2)\cr }\right]  =
   \left[\matrix{X_{n+1} & Y_{n+1} \cr Z_{n+1} & U_{n+1} \cr}\right]
   \left[\matrix{p_n(y_1)\cr p_n(y_2)\cr }\right] \;,
   \hskip 77pt
$$
whence the recurrence relations for the $X$'s, $Y$'s, $Z$'s, and
$U$'s:
$$\eqalign{
    a_{n+1}X_{n+1} &= y_1-b_n-a_nU_n/\delta_n\;, \cr
    a_{n+1}U_{n+1} &= y_2-b_n-a_nX_n/\delta_n\;, \cr
    a_{n+1}Y_{n+1} &= a_nY_n/\delta_n\;, \cr
    a_{n+1}Z_{n+1} &= a_nZ_n/\delta_n\;, \cr
    a_{n+1}^2\delta_{n+1} &= (y_1-b_n)(y_2-b_n)
    -a_n[(y_1-b_n)X_n+(y_2-b_n)U_n]/\delta_n+a_n^2/\delta_n\;. \cr
}$$
Let $\delta_n=\Theta_n/\Theta_{n-1}$ (yes, this is the $\Theta_n$
which will  appear in (9)), then,
with $\Theta_{n-1}X_n=\Upsilon_n$ and
     $\Theta_{n-1}U_n=\chi_n$,
$$\eqalign{
   a_{n+1}\Upsilon_{n+1} &= (y_1-b_n)\Theta_n-a_n\chi_n\;, \cr
   a_{n+1}\chi_{n+1} &= (y_2-b_n)\Theta_n-a_n\Upsilon_n\;, \cr
   a_{n+1}^2\Theta_{n+1} &= (y_1-b_n)(y_2-b_n)\Theta_n
    -a_n[(y_1-b_n)\Upsilon_n+(y_2-b_n)\chi_n]+a_n^2\Theta_{n-1}\;. \cr
}$$
Remark that $a_n\Theta_{n-1}Y_n$ and
            $a_n\Theta_{n-1}Z_n$ are independent of $n$.\hfill\break
These recurrence relations for $\Upsilon_n$, $\chi_n$ and $\Theta_n$
are     exactly the recurrence {\hbox{relations}} satisfied by
 {\it products\/}
of solutions of (2) at $y_1$ and $y_2$! Indeed, if
$a_{n+1}\xi_{n+1}=(y_1-b_n)\xi_n-a_n\xi_{n-1}$, and
$a_{n+1}\eta_{n+1}=(y_2-b_n)\eta_n-a_n\eta_{n-1}$, one finds for
$[\xi_n\eta_{n-1}, \xi_{n-1}\eta_n, \xi_n\eta_n]$ exactly the
recurrence for $[\Upsilon_n, \chi_n, \Theta_n]$ (Actually, this is
a recurrence of {\it fourth\/} order, one should work with vectors
$[\xi_n\eta_{n-1}, \xi_{n-1}\eta_n, \xi_n\eta_n, \xi_{n-1}\eta_{n-1}]$
and  $[\Upsilon_n, \chi_n, \Theta_n, \Theta_{n-1}]$).
 
\smallskip
Now, any solution of (2) is a combination of $p_n(x)$ and $q_n(x)$
defined by
$$ q_n(x)=\int_{{\rm Supp.\ }\mu}p_n(t)(x-t)^{-1}d\mu(t)=
      1/(\gamma_n x^{n+1})+\cdots\;,$$
a useful well-known identity is
$$ p_nq_{n-1}-p_{n-1}q_n=\gamma_n/\gamma_{n-1}=1/a_n,$$
so, $\xi_n$ is some combination of $p_n(y_1)$ and $q_n(y_1)$,
    $\eta_n$ is some combination of $p_n(y_2)$ and $q_n(y_2)$, and
$$\hskip -18pt\matrix{
\Upsilon_n &=&\alpha p_n(y_1)p_{n-1}(y_2) &+& \beta p_n(y_1)q_{n-1}(y_2)
  &+& \gamma q_n(y_1)p_{n-1}(y_2) &+&\delta q_n(y_1)q_{n-1}(y_2)\;, \cr
\chi_n &=&\alpha p_{n-1}(y_1)p_n(y_2) &+& \beta p_{n-1}(y_1)q_n(y_2)
  &+& \gamma q_{n-1}(y_1)p_n(y_2) &+&\delta q_{n-1}(y_1)q_n(y_2)\;, \cr
\Theta_n &=&\alpha p_n(y_1)p_n(y_2) &+& \beta p_n(y_1)q_n(y_2)
  &+& \gamma q_n(y_1)p_n(y_2) &+&\delta q_n(y_1)q_n(y_2)\;. \cr
}$$
As we want the degrees of the left sides to remain bounded when $n$
increases, we must have $\alpha=0$. Let us show that $\delta=0$ as
well: from
$ p_n(y_1)=X_np_{n-1}(y_1)+Y_np_{n-1}(y_2)\;$,
$ \Theta_{n-1}p_n(y_1)-\Theta_{n-1}X_np_{n-1}(y_1)=$
$ \Theta_{n-1}p_n(y_1)-\Upsilon_n p_{n-1}(y_1)=$
$[\beta p_{n-1}(y_1)q_{n-1}(y_2)
  + \gamma q_{n-1}(y_1)p_{n-1}(y_2) +\delta q_{n-1}(y_1)q_{n-1}(y_2)]
p_n(y_1)-$
$[\beta p_n(y_1)q_{n-1}(y_2)
  + \gamma q_n(y_1)p_{n-1}(y_2) +\delta q_n(y_1)q_{n-1}(y_2)]
p_{n-1}(y_1)=$
$\gamma p_{n-1}(y_2)[q_{n-1}(y_1)p_n(y_1)-q_n(y_1)p_{n-1}(y_1)]+$
\hfill\break
$\delta q_{n-1}(y_2)[q_{n-1}(y_1)p_n(y_1)-q_n(y_1)p_{n-1}(y_1)]=$
$[\gamma p_{n-1}(y_2)+\delta q_{n-1}(y_2)]/a_n=$
$\Theta_{n-1}Y_np_{n-1}(y_2)\;$, or \hfill\break
$\gamma p_{n-1}(y_2)+\delta q_{n-1}(y_2)=$
$a_n\Theta_{n-1}Y_n p_{n-1}(y_2)=$ constant $p_{n-1}(y_2)\;$,
possible only if the constant (with respect to $n$)
$a_n\Theta_{n-1}Y_n =\gamma$, and $\delta=0$.
 
One finds similarly $a_n\Theta_{n-1}Z_n=\beta$.
 
Finally, at $n=0$,
$\Theta_0=\beta p_0(y_1)q_0(y_2)
      +   \gamma q_0(y_1)p_0(y_2) =
[\beta S(y_2)+\gamma S(y_1)]/\mu_0$,
and this relation between $S(y_1)$ and $S(y_2)$ is exactly (8)!
One has $U=\Theta_0$, $W=(y_2-y_1)(\beta-\gamma)/(2\mu_0)$, and
$V=-(\beta+\gamma)/(2\mu_0)$.

\hfill\QED
\bigskip
 
{\bf 4. Semi-classical measures  on $snul$s.}

\bigskip
Let us consider meromorphic Stieltjes functions $S$, corresponding
therefore to discrete (atomic [8]) measures.
From (8), or  from the  equivalent form
$\beta S(y_2)+\gamma S(y_1)=\mu_0 U$, we have a recurrence
$$\beta(x_k,y_k)S(y_{k+1})+\gamma(x_k,y_k)S(y_k)=\mu_0 U(x_k),$$
showing that poles occur at some available lattice $\ldots, y_k,
y_{k+1},\ldots$ with residues (masses of the measure) satisfying
$$\beta(x_k,y_k)\mu(y_{k+1})+\gamma(x_k,y_k)\mu(y_k)=0,\eqno(12)$$
a Pearson-like equation like this is discussed in [29,30].
 
The masses will usually not make an infinite convergent sequence,
so that they must remain in finite number if one wants a discrete
measure. This is possible only if some value of $\beta(x,y)$ on a
lattice, say $\beta(x_0,y_0)$ vanishes, so that we may start the
nonzero masses at $y_1$ (and put $\mu(y_0)=0$), and also if some
value of $\gamma(x,y)$ vanishes on the same lattice, say
$\gamma(x_N,y_N)=0$, so that we stop the nonzero masses at $y_N$
(and put $\mu(y_{N+1})=0$), see [29] (the ``uninteresting case'', p.655).
 
For instance, if $V=0$ in (8), we simply have equal masses:
$\beta(x,y)=-\gamma(x,y)=W(x)$, and $W$ must vanish at two points of
the $x-$lattice at least. If this does not happen, the measure is
approximated by a discrete measure with many
small equal masses at more and more
lattice points, and tends towards the limit distribution of these
lattice points. Take for instance (6) with $x_c=y_c=0$,
$\zeta\sqrt{A}=
 \zeta\sqrt{C}=1/2$, and $q=\exp(i\theta)$, then, $x_k=\cos k\theta$,
$y_k=\cos(k-1/2)\theta$, the lattice VI of [25], distributed (if
$\theta/\pi$ is not rational) like $d\mu(x)=(1-x^2)^{-1/2}dx$: we
recover the Chebyshev polynomials! This is not surprising, as these
polynomials have interesting properties with respect to the
$\D-$operator [6, 12].
\bigskip
 
\vfill\eject
 
{\bf 5. Dense discrete measures  and combinations of Chebyshev
polynomials.}

\bigskip
We keep  $x_k=\cos k\theta$,
$y_k=\cos(k-1/2)\theta$, the lattice VI of [25] with
$\theta/\pi$  not rational, and try to find a discrete measure with
jumps at each of these $y_k$.
 ``Strange''
supports, er, {\it carriers\/},
 have been well worked [11,20,23,31,32,33], there is a reference to
Stieltjes himself in p.202 of [18].
 
As semi-classical orthogonal polynomials do not seem to be related
to dense discrete measures, one tries combinations of the simplest
such items, i.e., Chebyshev polynomials:
 
\medskip
 
\noindent{\bf Theorem.\/} {\sl Let the measure $\mu$ be discrete
with jumps $\mu(y_k)=1/(k-1/2)^2$ at $y_k=\cos((k-1/2)\theta)$,
$k=\ldots,-2,-1,0,1,2,\ldots$, with $\theta/\pi$ irrational.
Then, the orthonormal polynomials are
$$ p_0=\pi^{-1},\hskip 355pt $$
$$\displaylines{
 \hskip -8pt
  p_n = \dfrac
  {(\varepsilon_{2n-1}+\iota_{2n-1})T_n
   -(-1)^{\lfloor\xi_{2n-1}\theta/(2\pi)\rfloor}
  \iota_{2n-1} T_{|n-\xi_{2n-1}|}
    +(-1)^{\lfloor\eta_{2n-1}\theta/(2\pi)\rfloor}\varepsilon_{2n-1}
     T_{|n-\eta_{2n-1}|}
   }
   {\sqrt{
   2\pi^2 \varepsilon_{2n} \iota_{2n} (\varepsilon_{2n-1}+\iota_{2n-1})}}
\;,\cr
\hfill \ n\ge 1 \qquad (13)    \cr
}
$$
where $\xi_j$ is the value of $p$ which minimizes
 $p\theta/(2\pi)-\lfloor p\theta/(2\pi)\rfloor$ on $p=1,2,\ldots,j$,
\hfill\break
 $\eta_j$ is the value of $p$ which minimizes
 $\lfloor p\theta/(2\pi)\rfloor+1-p\theta/(2\pi)$ on $p=1,2,\ldots,j$,
\hfill\break
$($ where $\lfloor x \rfloor$ is the largest integer smaller or equal
 than $x)$, $\varepsilon_j=\xi_j\theta/(2\pi)-
                 \lfloor   \xi_j\theta/(2\pi)\rfloor$,
$ \iota_j =  \lfloor   \eta_j\theta/(2\pi)\rfloor  +1-
                       \eta_j\theta/(2\pi)$.
 
}
 
\medskip
 
Remark that $\varepsilon_j$ and $\iota_j$ are positive and decreasing
with $j$. These $\xi$'s and $\eta$'s are known as denominators of
remarkable rational approximants to the irrational number $\theta/(2\pi)$
(``Nebenn\"aherungsbr\"uche'' in [26] \S~16) and are linked to the continued
fraction expansion of $\theta/(2\pi)$. Each new $\xi$ or $\eta$ is the
sum of the two last ones: if $\varepsilon_j>\iota_j$,
$\ell=\xi_\ell=\xi_j+\eta_j$ and $\varepsilon_\ell=\varepsilon_j-
\iota_j$, $\iota_\ell=\iota_j$;
                          if $\varepsilon_j<\iota_j$,
$\ell=\eta_\ell=\xi_j+\eta_j$ and $\iota_\ell=\iota_j-\varepsilon_j$,
$\varepsilon_\ell=\varepsilon_j$.
 For instance, with the golden ratio
$\theta/(2\pi)=(5^{1/2}+1)/2 =
 1.6180339887498948481\ldots$, one encounters the\dots Fibonacci
numbers:
 
% GenTeX page 54
$$\hskip 31pt \vbox{
   \offinterlineskip
   \halign{
   \strutt \vrule \hfil\ $#$\  & \vrule \hfil\ $#$\ \hfil & \vrule
    \hfil\ $#$\ \hfil  \vrule \cr
   \noalign{\hrule}
   \xi  {\hbox{\rm \ or\ }} \eta\ \ & \varepsilon=\xi\theta/(2\pi)-
                       \lfloor\xi\theta/(2\pi)\rfloor &
      \iota=   \lfloor\eta\theta/(2\pi)\rfloor +1-
                      \eta\theta/(2\pi) \cr
   \noalign{\hrule}
 1 &  0.6180339887498948481 &   0.3819660112501051518  \cr
 2 &  0.2360679774997896963 &                           \cr
 3 &                        & 0.1458980337503154555     \cr
 5 &  0.0901699437494742407 &                           \cr
 8 &                        & 0.0557280900008412147     \cr
 13 &  0.0344418537486330259&                           \cr
 21 &                       & 0.0212862362522081887     \cr
 34 &  0.0131556174964248372&                           \cr
 55 &                       & 0.0081306187557833515     \cr
   \noalign{\hrule}
            }
\medskip
\centerline{}
 
       } $$
 
With $\theta/(2\pi)=\sqrt{2}$, one finds new values of $\xi$ at
1, 3, 5, 17, 29,\dots, new $\eta$ values at 1, 2, 7, 12, 41,\dots
 
Remark also that the present discrete measure is definitely {\it not\/}
semi-classical, as $\mu(y_{k+1})/\mu(y_k) =
((k-1/2)/(k+1/2))^2$ is not a rational function of
$x_k=\cos k\theta$ and $y_k=\cos((k-1/2)\theta)$, as it should have been
according to (12).

\smallskip
\noindent {\bf Proof of the theorem.\/} In order to show that a form
like $p_n=A_n T_n+B_nT_{n-\xi_{2n-1}}+C_n T_{n-\eta_{2n-1}}$ is valid
 (using
$T_{|p|}=T_p$), one must show that the scalar product
$$\displaylines{
  (p_n,T_m) =\sum_{k=-\infty}^\infty \mu(y_k) p_n(y_k)T_m(y_k) \hfill\cr
 \hskip 39pt =\dfrac12\sum_{k=-\infty}^\infty \mu(y_k)\left[
   A_n(T_{n+m}(y_k)+T_{n-m}(y_k))
 +B_n(T_{n+m-\xi_{2n-1}}(y_k)+T_{n-m-\xi_{2n-1}}(y_k)) \right. \cr
 \hfill \left.
 +C_n(T_{n+m-\eta_{2n-1}}(y_k)+T_{n-m-\eta_{2n-1}}(y_k))\right] \cr
        } $$
vanishes for $m=0,1,\ldots, 2n-1$, and has the value
$2^{n-1}/\gamma_n=1/A_n$ when $m=2n$.
 
Let
$$ \eqalign{
  \tau_p &= \sum_{k=-\infty}^\infty \mu(y_k)T_p(y_k) \cr
         &= \sum_{k=-\infty}^\infty (k-1/2)^{-2} \cos(p(k-1/2)\theta) \cr
        &= 8[ \cos(p\theta/2) +\cos(3p\theta/2)/9 +\cos(5p\theta/2)/25+\cdots]
          \cr
       &= 2\pi^2 (-1)^{\lfloor p\theta/(2\pi)\rfloor}
            \left( \lfloor p\theta/(2\pi)\rfloor +1/2-p\theta/(2\pi)
           \right) \hskip 105pt (14')\cr
       &= -2\pi^2 (-1)^{\lceil p\theta/(2\pi)\rceil}
            \left( \lceil p\theta/(2\pi)\rceil -1/2-p\theta/(2\pi)
           \right)\;, \hskip 95pt (14'{}')\cr
             }
$$
from elementary Fourier series,
where $\lceil x \rceil$ is the smallest integer larger or equal than $x$.
The form of $(14'{}')$ will sometimes be more convenient than $(14')$.
The importance of using Chebyshev moments of measures with discrete
masses has been shown by Pr\'evost ([27,28]), remark in particular
that $\tau_p$ does {\it not\/} tend to zero when $p\to\infty$, as it
should with absolutely continuous measures.
So, we have $(p_n,T_m)=(A_n\tau_{n+m}+B_n\tau_{n+m-\xi_{2n-1}}+
C_n\tau_{n+m-\eta_{2n-1}}
       +A_n\tau_{n-m}+B_n\tau_{n-m-\xi_{2n-1}}+C_n\tau_{n-m-\eta_{2n-1}})/2$.
Remark that $\tau_0=\pi^2$, so that $p_0=(\tau_0)^{-1/2}=\pi^{-1}$.
For $n\ge 1$,
let us show that it is possible to find $A_n$, $B_n$, and $C_n$ such
that $A_n\tau_N+B_n\tau_{N-\xi_{2n-1}}+C_n\tau_{N-\eta_{2n-1}}=0$ for
$N=1,2,\ldots, 2n-1$:\hfill\break
let $\rho=\theta/(2\pi)$.
 As
$0<\xi_{2n-1}\rho-\lfloor\xi_{2n-1}\rho\rfloor \le
       N\rho-\lfloor    N\rho\rfloor$,
$(N-\xi_{2n-1})\rho-1 < \lfloor N\rho\rfloor
  -\xi_{2n-1}\rho <
    \lfloor N\rho\rfloor
  -\lfloor\xi_{2n-1}\rho\rfloor \le
 (N-\xi_{2n-1})\rho$, one has
$   \lfloor N\rho\rfloor
  -\lfloor\xi_{2n-1}\rho\rfloor =
   \lfloor (N-\xi_{2n-1})\rho\rfloor$.
                      As
$\lfloor\eta_{2n-1}\rho\rfloor+1-\eta_{2n-1}\rho \le
   N\rho\rfloor +1 -     N\rho <1$,
$(N-\eta_{2n-1})\rho   \le \lfloor N\rho\rfloor
  -\lfloor\eta_{2n-1}\rho\rfloor <
  -\lfloor\eta_{2n-1}\rho\rfloor
  +N\rho <
    1+
 (N-\eta_{2n-1})\rho$, one has
$   \lfloor N\rho\rfloor
  -\lfloor\eta_{2n-1}\rho\rfloor =
   \lceil  (N-\eta_{2n-1})\rho\rceil $.
 
\smallskip\noindent
In summary:
$$  \lfloor N\rho\rfloor   -\lfloor\xi_{2n-1}\rho\rfloor =
   \lfloor (N-\xi_{2n-1})\rho\rfloor\;,\;
    \lfloor N\rho\rfloor  -\lfloor\eta_{2n-1}\rho\rfloor =
   \lceil  (N-\eta_{2n-1})\rho\rceil,  \eqno(15) $$
for $N=1,2,\ldots, 2n-1$.
 
So,
$\tau_{N-\xi_{2n-1}} = 2\pi^2 (-1)^{\lfloor N\rho\rfloor -
                               \lfloor \xi_{2n-1}\rho\rfloor }
  \left(
                               \lfloor N\rho\rfloor -
                               \lfloor \xi_{2n-1}\rho\rfloor
     +1/2
          - N\rho\right.$ $\left. +  \xi_{2n-1}\rho\right) =
  (-1)^{\lfloor \xi_{2n-1}\rho\rfloor } \tau_N
  +2\pi^2
                         (-1)^{\lfloor N\rho\rfloor -
                               \lfloor \xi_{2n-1}\rho\rfloor }
            \varepsilon_{2n-1}$,
\hfill\break
$\tau_{N-\eta_{2n-1}} = -2\pi^2 (-1)^{\lfloor N\rho\rfloor -
                               \lfloor \eta_{2n-1}\rho\rfloor }
  \left(
                               \lfloor N\rho\rfloor -
                               \lfloor \eta_{2n-1}\rho\rfloor
     -1/2
          - N\rho\right.$ $\left. +  \eta_{2n-1}\rho\right) =
\hfill\break
 -(-1)^{\lfloor \eta_{2n-1}\rho\rfloor } \tau_N
  +2\pi^2
                         (-1)^{\lfloor N\rho\rfloor -
                               \lfloor \xi_{2n-1}\rho\rfloor }
            \iota_{2n-1}$,
\hfill\break
 $A_n\tau_N+B_n\tau_{N-\xi_{2n-1}}+C_n\tau_{N-\eta_{2n-1}}=
[A_n
 +(-1)^{\lfloor \xi_{2n-1}\rho\rfloor } B_n
 -(-1)^{\lfloor \eta_{2n-1}\rho\rfloor } C_n]\tau_N
  +2\pi^2(-1)^{\lfloor N\rho\rfloor}
    [(-1)^{\lfloor \xi_{2n-1}\rho\rfloor }\varepsilon_{2n-1} B_n +
     (-1)^{\lfloor \eta_{2n-1}\rho\rfloor }\iota_{2n-1} C_n]$,
which vanishes indeed for $N=n\pm m=1,2,\ldots, 2n-1$ if
$B_n=-K_n (-1)^{\lfloor \xi_{2n-1}\rho\rfloor }\iota_{2n-1}$,
$C_n =K_n(-1)^{\lfloor \eta_{2n-1}\rho\rfloor }\varepsilon_{2n-1}$,
and $A_n=K_n(\varepsilon_{2n-1}+\iota_{2n-1})$.
 
\smallskip\noindent
We now have to look at
$(p_n,T_n)=(A_n\tau_{2n}+B_n\tau_{2n-\xi_{2n-1}}+C_n\tau_{2n-\eta_{2n-1}}
       +A_n\tau_0+B_n\tau_{\xi_{2n-1}}+C_n\tau_{\eta_{2n-1}})/2=1/A_n$.
 
\item{1.} If $\xi_{2n}=\xi_{2n-1}<2n$ and $\eta_{2n}=\eta_{2n-1}<2n$, nothing
changes in the {\hbox{evaluation}} of
 $A_n\tau_N+ B_n\tau_{N-\xi_{2n-1}}+ C_n\tau_{N-\eta_{2n-1}}$ when we
replace $N$ by $2n$, so
 $A_n\tau_{2n}+ B_n\tau_{2n-\xi_{2n-1}}+ C_n\tau_{2n-\eta_{2n-1}}=0$, and
we have only to look at
  $A_n\tau_0+ B_n\tau_{\xi_{2n-1}}+ C_n\tau_{\eta_{2n-1}} $.
One has $\tau_0=\pi^2$,
 $\tau_{\xi_{2n-1}}=2\pi^2(-1)^{\lfloor \xi_{2n-1}\rho\rfloor}
   (1/2-\varepsilon_{2n-1})$,
 $\tau_{\eta_{2n-1}}=-2\pi^2(-1)^{\lfloor \eta_{2n-1}\rho\rfloor}
   (1/2-\iota_{2n-1})$,
yielding
  $A_n\tau_0+ B_n\tau_{\xi_{2n-1}}+ C_n\tau_{\eta_{2n-1}} =
  4K_n\pi^2\varepsilon_{2n-1}\iota_{2n-1}$, which must be equal to
 $2/A_n=2/(K_n(\varepsilon_{2n-1}+\iota_{2n-1}))$, whence
$K_n = [(\varepsilon_{2n-1}+\iota_{2n-1})/(2 \pi^2
         \varepsilon_{2n-1} \iota_{2n-1})  ]^{1/2}$,
and this gives  (13), as
we still have $\varepsilon_{2n}=\varepsilon_{2n-1}$ and
$\iota_{2n}=\iota_{2n-1}$.
 
\item{2.} $\xi_{2n}=2n$ or $\eta_{2n}=2n$, which happens only if
$\xi_{2n-1}+\eta_{2n-1} = 2n$.
Then,
  $A_n\tau_0+ B_n\tau_{\xi_{2n-1}}+ C_n\tau_{\eta_{2n-1}} +
  A_n\tau_{2n}+ B_n\tau_{2n-\xi_{2n-1}}+ C_n\tau_{2n-\eta_{2n-1}} =
A_n(\tau_0+\tau_{2n})+ (B_n+C_n)(\tau_{\xi_{2n-1}}+ \tau_{\eta_{2n-1}})$.
\hfill\break
An interesting consequence of (15) is that $\lfloor \xi_{2n-1}\rho\rfloor$
and
  $\lfloor \eta_{2n-1}\rho\rfloor$ have now the same evenness: if we
subtract the two equations of (15) with $N=n=(\xi_{2n-1}+\eta_{2n-1})/2$,
one finds
$-\lfloor \xi_{2n-1}\rho\rfloor+ \lfloor \eta_{2n-1}\rho\rfloor=
  \lfloor (\eta_{2n-1}-\xi_{2n-1})\rho/2\rfloor-
  \lceil (\xi_{2n-1}-\eta_{2n-1})\rho/2\rceil$,
which is an even integer, as $\lfloor x\rfloor= -\lceil -x\rceil$
([15] \S~1.2.4).
\hfill\break
So, let $\sigma= (-1)^{\lfloor\xi_{2n-1}\rho\rfloor}=
                 (-1)^{\lfloor\eta_{2n-1}\rho\rfloor}$. One has
$B_n+ C_n= \sigma K_n (\varepsilon_{2n-1}- \iota_{2n-1})$,
$\tau_{\xi_{2n-1}}+ \tau_{\eta_{2n-1}}= 2\pi^2 \sigma(1/2-\varepsilon_{2n-1})-
2\pi^2 \sigma (1/2-\iota_{2n-1})= 2\pi^2\sigma
(\iota_{2n-1}-\varepsilon_{2n-1})$.
 
\itemitem{2a.} If $\xi_{2n}=2n$, $\varepsilon_{2n}=
                                  \varepsilon_{2n-1}- \iota_{2n-1}$,
$2n\rho- \lfloor 2n\rho\rfloor=
 \xi_{2n-1}\rho- \lfloor\xi_{2n-1}\rho\rfloor-
 (\lfloor\eta_{2n-1}\rho\rfloor+1- \eta_{2n-1}\rho)$, so,
$\lfloor 2n\rho\rfloor=
 \lfloor\xi_{2n-1}\rho\rfloor+
  \lfloor\eta_{2n-1}\rho\rfloor-1$,
$\tau_{2n}= 2\pi^2 (-1)^{\lfloor 2n\rho\rfloor} (1/2-\varepsilon_{2n})=
-2\pi^2(1/2-\varepsilon_{2n})$,
$A_n(\tau_0+\tau_{2n})+ (B_n+C_n)(\tau_{\xi_{2n-1}}+ \tau_{\eta_{2n-1}})=
4\pi^2 K_n \varepsilon_{2n} \iota_{2n-1}$, whence (13), as one still
has $\iota_{2n}=\iota_{2n-1}$.
 
\itemitem{2b.} If $\eta_{2n}=2n$, $\iota_{2n}=
                                  \iota_{2n-1}-\varepsilon_{2n-1}$,
$\lfloor 2n\rho\rfloor+1- 2n\rho=
  \lfloor\eta_{2n-1}\rho\rfloor+1- \eta_{2n-1}\rho-
(\xi_{2n-1}\rho- \lfloor\xi_{2n-1}\rho\rfloor)$, so
$\lfloor 2n\rho\rfloor=
  \lfloor\eta_{2n-1}\rho\rfloor+
 \lfloor\xi_{2n-1}\rho\rfloor$,
$\tau_{2n}= 2\pi^2 (-1)^{\lfloor 2n\rho\rfloor} (-1/2+\iota_{2n})=
-2\pi^2(1/2-\iota_{2n})$,
$A_n(\tau_0+\tau_{2n})+ (B_n+C_n)(\tau_{\xi_{2n-1}}+ \tau_{\eta_{2n-1}})=
4\pi^2 K_n \iota_{2n} \varepsilon_{2n-1}$, whence (13), as one still
has $\varepsilon_{2n}=\varepsilon_{2n-1}$.
\hfill\QED
\bigskip

A numerical check has been performed with $\rho=$
$\theta/(2\pi)=(5^{1/2}+1)/2 =
 1.6180339887498948481\ldots$, the recurrence coefficients have been
computed by Gautschi's {\tt sti} (Stieltjes, of course!) subroutine [9].
 
For each $n$, one  compares the computed $\pi \gamma_n 2^{-n} =
1/(2^n a_1\ldots a_n)$ with the formula predicted from (13), i.e.,
$((\varepsilon_{2n-1}+\iota_{2n-1})/(8\varepsilon_{2n}\iota_{2n}))^{1/2}$.
The agreement is satisfactory, taking into account that all the series
involving $\mu(y_k)=1/(k-1/2)^2$ have been truncated to 20000 terms,
so that relative errors of about $10^{-3}$ may be expected.

We have longer and longer intervals where $a_n=1/2$ and $b_n=0$. On
this example, $\liminf_{n\to\infty} a_n >0$, although the vanishing of
this $\liminf$ could have been expected from pure discrete (atomic)
measures [8], but other results have been published about singular
measures [7,11,13,14,18,20,23,31,32,33].

\bigskip\bigskip

\noindent{\bf Acknowledgements.}
\medskip
Many thanks to P.\ Barrucand, P.\ Bulens, T.S.\ Chihara, J.\ Dombrowski,
 M.\ Ismail, J.\ Meinguet,
 P.\ Nevai, A.\ Ronveaux, W.\ Van Assche,
 for kind words and information.
 
Many thanks to W.A.\ Al-Salam too, who organizes a preprint
repository (several preprints given in the references list come from
there) at the anonymous ftp site {\tt euler.math.ualberta.ca}.
 
The saying of William of Orange, the most remarkable statesman
of the $16^{\hbox{\sevenrm{th}}}$ century [3] (could be compared
to N.\ Mandela nowadays), applies quite well to scientific research
(try to explain that to a contemporary state(?)sman).
 
Many thanks to the Organizing Committees of the meeting
{\bf TJS94\/} too!
 
\vfill\eject

% GenTeX page 54
$$\hskip 1pt \vbox{
   \offinterlineskip
   \halign{
   \strutt \vrule \hfil\ $#$\  & \vrule \hfil\ $#$\  &
  \vrule \hfil\ $#$\  & \vrule \hfil\ $#$\ & \vrule \hfil\ $#$\ &
  \vrule \hfil\ $#$\  & \vrule \hfil\ $#$\ & \vrule \hfil $#$ \hfil &
    \vrule \hfil\ $#$\hfil   \vrule \cr
   \noalign{\hrule}
   n & a_n & b_n &
   \varepsilon_{2n-1} & \iota_{2n-1} &
   \varepsilon_{2n} & \iota_{2n} &
  \dfrac1{2^n a_1\ldots a_n} &
  \sqrt{\dfrac{\varepsilon_{2n-1}+\iota_{2n-1}}
              {8\varepsilon_{2n}\iota_{2n}}} \cr
   \noalign{\hrule}
 0 &  & 0.2361 & & & & & & \cr
   1& 0.4247&-0.5451& 0.6180& 0.3820& 0.2361& 0.3820&    1.177&    1.177\cr
   2& 0.5\ \ \ & 0.6180& 0.2361& 0.1459& 0.2361& 0.1459&    1.177&    1.177\cr
   3& 0.3931&-0.3090& 0.0902& 0.1459& 0.0902& 0.1459&    1.498&    1.498\cr
   4& 0.3090& 0\ \ \ \ \ & 0.0902& 0.1459& 0.0902& 0.0557&    2.423&    2.423\cr
   5& 0.6360& 0\ \ \ \ \ & 0.0902& 0.0557& 0.0902& 0.0557&    1.905&    1.905\cr
   6& 0.5\ \ \ &-0.3090& 0.0902& 0.0557& 0.0902& 0.0557&    1.905&    1.905\cr
   7& 0.3931& 0.3090& 0.0344& 0.0557& 0.0344& 0.0557&    2.423&    2.423\cr
   8& 0.5\ \ \ & 0\ \ \ \ \ & 0.0344& 0.0557& 0.0344& 0.0557&    2.424&    2.423\cr
   9& 0.5\ \ \ & 0\ \ \ \ \ & 0.0344& 0.0557& 0.0344& 0.0557&    2.423&    2.423\cr
  10& 0.5\ \ \ & 0.3090& 0.0344& 0.0557& 0.0344& 0.0557&    2.423&    2.423\cr
  11& 0.3931&-0.3090& 0.0344& 0.0213& 0.0344& 0.0213&    3.083&    3.082\cr
  12& 0.5\ \ \ & 0\ \ \ \ \ & 0.0344& 0.0213& 0.0344& 0.0213&    3.083&    3.082\cr
    & \cdots& \cdots& \cdots& \cdots& \cdots& \cdots&  \cdots & \cdots\cr
  17& 0.3090& 0\ \ \ \ \ & 0.0344& 0.0213& 0.0132& 0.0213&    4.989&    4.988\cr
  18& 0.6360& 0\ \ \ \ \ & 0.0132& 0.0213& 0.0132& 0.0213&    3.922&    3.921\cr
  19& 0.5\ \ \ & 0\ \ \ \ \ & 0.0132& 0.0213& 0.0132& 0.0213&    3.922&    3.921\cr
    & \cdots& \cdots& \cdots& \cdots& \cdots& \cdots&  \cdots & \cdots\cr
  27& 0.5\ \ \ &-0.3091& 0.0132& 0.0213& 0.0132& 0.0213&    3.922&    3.921\cr
  28& 0.3930& 0.3091& 0.0132& 0.0081& 0.0132& 0.0081&    4.989&    4.988\cr
  29& 0.5\ \ \ & 0\ \ \ \ \ & 0.0132& 0.0081& 0.0132& 0.0081&    4.989&    4.988\cr
    & \cdots& \cdots& \cdots& \cdots& \cdots& \cdots&  \cdots & \cdots\cr
  44& 0.5\ \ \ & 0.3091& 0.0132& 0.0081& 0.0132& 0.0081&    4.989&    4.988\cr
  45& 0.3930&-0.3091& 0.0050& 0.0081& 0.0050& 0.0081&    6.348&    6.344\cr
  46& 0.5\ \ \ & 0\ \ \ \ \ & 0.0050& 0.0081& 0.0050& 0.0081&    6.348&    6.344\cr
    & \cdots& \cdots& \cdots& \cdots& \cdots& \cdots&  \cdots & \cdots\cr
  71& 0.5\ \ \ & 0\ \ \ \ \ & 0.0050& 0.0081& 0.0050& 0.0081&    6.348&    6.344\cr
  72& 0.3089& 0\ \ \ \ \ & 0.0050& 0.0081& 0.0050& 0.0031&   10.277&   10.265\cr
  73& 0.6361& 0\ \ \ \ \ & 0.0050& 0.0031& 0.0050& 0.0031&    8.078&    8.070\cr
  74& 0.5\ \ \ & 0\ \ \ \ \ & 0.0050& 0.0031& 0.0050& 0.0031&    8.078&    8.070\cr
   \noalign{\hrule}
            }
\medskip
\centerline{}
 
       } $$
 
\vfill\eject
 
\noindent{\bf References.}
 
\frenchspacing      \parindent=20pt
 
\medskip

\item{[1]} W.A.\ Al-Salam, Characterization theorems for orthogonal
    polynomials, pp. 1-24 {\it in\/} {\sl Orthogonal Polynomials:
    Theory and Practice\/}, (P.\ Nevai, editor) {\sl NATO ASI Series
   C: Math.\ and Phys.\ Sciences\/} {\bf 294\/}, Kluwer, Dordrecht, 1990.
 
\item{[2]} R.~Askey, J.~Wilson, Some basic hypergeometric
      orthogonal polynomials that generalize Jacobi polynomials,
      {\sl Mem.\ A.\ M.\ S.\/} {\bf 54\/} no.\ 319, 1985.
 
\item{[3]} R.~Avermaete, {\sl Guillaume d'Orange, dit le Taciturne,
      1533-1584.} Payot, Paris, 1939, reprinted 1984.
 
\item{[4]}  S.~Bonan, P.~Nevai, Orthogonal polynomials and their
       derivatives,I, {\sl J. Approx. Theory\/} {\bf 40\/} (1984),
       134-147.
 
\item{[5]}  S.~Bonan, D.S.~Lubinsky, P.~Nevai, Orthogonal polynomials
 and their
       derivatives,II, {\sl SIAM J. Math. An.\/} {\bf 18\/} (1987),
       1163-1176.

\item{[6]} B.M.~Brown, M.E.H.~Ismail, A right inverse of the
     Askey-Wilson operator. Preprint.
 
\item{[7]} R.~Del Rio, N.~Makarov, B.~Simon, Operators with singular
     continuous spectrum:~II, rank one operators,
    {\sl Comm.\ Math.\ Phys.\/} {\bf 165\/} (1994) 59-67.
 
\item{[8]} J.\ Dombrowski, Tridiagonal matrix representations of
    cyclic self-adjoint operators, {\sl Pacific J.\ Math.\ }
    {\bf 114\/} (1984) 325-334. II, {\it ibidem\/} {\bf 120\/}
    (1985) 47-53.
 
\item{[9]}  W.\ Gautschi, Algorithm 726.\
    ORTHPOL: a package of routines for generating orthogonal polynomials
    and Gauss-type quadrature rules, {\sl ACM Trans.\ Math.\ Soft.\ }
    {\bf 20\/} (1994) 21-62.
 
\item{[10]}  E.~Hendriksen, H.~van Rossum, Semi-classical orthogonal
     polynomials,
                                       pp. 354-361 {\it in\/}
      {\sl Polyn\^omes Orthogonaux et Applications, Proceedings,
      Bar-le-Duc 1984\/}, (C.Brezinski \& al., editors),
      {\sl Lecture Notes Math.\/} {\bf 1171\/}, Springer, Berlin 1985.
 
\item{[11]} A.~Iserles, From Schr\"odinger spectra to orthogonal
      polynomials via a functional equation, preprint DAMTP Cambridge.
 
\item{[12]}
     M.E.H.~Ismail, R.~Zhang, Diagonalization
     of certain integral operators. M.E.H.~Ismail, M.~Rahman, R.~Zhang,
     Diagonalization of certain integral operators II. Preprints.
 
\item{[13]} S.Ya.~Jitomirskaya, Anderson localization for the almost
     Mathieu equation: a nonperturbative proof,
    {\sl Comm.\ Math.\ Phys.\/} {\bf 165\/} (1994) 49-57.
 
\item{[14]} S.~Jitomirskaya, B.~Simon,
                                              Operators with singular
     continuous spectrum:~III, almost periodic Schr\"odinger operators,
    {\sl Comm.\ Math.\ Phys.\/} {\bf 165\/} (1994) 201-205.
 
\item{[15]} D.E.~Knuth, {\sl The Art of Computer Programming\/} {\bf 1\/}:
     {\sl Fundamental Algorithms.} Addison-Wesley, Reading, 1968.
 
\item{[16]}  R.~Koekoek, R.F.~Swarttouw, The Askey-scheme of hypergeometric
      polynomials and its $q-$analogue,
      preprint T.U.\ Delft.
 
\item{[17]} T.H.~Koornwinder,
  {\sl Compact quantum groups and $q$-special functions},
  in {\sl Representations of Lie groups and quantum groups},
  V. Baldoni \& M. A. Picardello (eds.),
  Pitman Research Notes in Mathematics Series 311,
  Longman Scientific \& Technical, 1994.

\item{[18]} H.~Kunz, B.~Souillard, Sur le spectre des op\'erateurs
   aux diff\'erences finies al\'eatoires, {\sl Commun.\ Math.\ Phys.\/}
   {\bf 78\/} (1980) 201-246.

\item{[19]}  E.~Laguerre, Sur la r\'eduction en fractions continues d'une
      fraction qui satisfait \`a une \'equation diff\'erentielle lin\'eaire
      du premier ordre dont les coefficients sont rationnels,
      {\sl J. Math. Pures Appl. (4)\/} {\bf 1\/} (1885), 135-165 =
      pp. 685-711 {\it in\/}
      {\sl Oeuvres\/}, Vol.II, Chelsea, New-York 1972.
 
\item{[20]} D.S.~Lubinsky, Jump distributions on $[-1,1]$ whose orthogonal
     polynomials have leading coefficients with given asymptotic behaviour,
     {\sl Proc.\ A.\ M.\ S.\/} {\bf 104\/} (1988) 516-524.
 
\item{[21]} A.P.~Magnus, Associated Askey-Wilson polynomials as
   Laguerre-Hahn orthogonal polynomials, pp.~261-278 {\it in\/}
   M.~Alfaro {\it et al.\/}, editors: {\sl Orthogonal Polynomials and
   their Applications, Proceedings, Segovia 1986. Springer Lecture
   Notes Math.\/} {\bf 1329\/}, Springer, Berlin, 1988.
 
\item{[22]} P.\ Maroni, Une caract\'erisation des  polyn\^omes
  orthogonaux semi-classiques, {\sl C.R.\ Acad.\ Sci.\ Paris ser.\ \/}
   1 {\bf 301\/} (1985) 269-272.
 
\item{[23]} S.N.\ Naboko, S.I.\ Yakovlev, The discrete Schrodinger
     operator. The point spectrum lying on the continuous spectrum,
     {\sl Algebra i Analys\/} {\bf 4\/} (1992) 183-195 (in Russian) =
  {\sl St.\ Petersburg Math.\ Journal\/} {\bf 4\/} (1993) 559-568.
 
\item{[24]} A.F.\ Nikiforov, V.B.\ Uvarov, {\sl Classical orthogonal
  polynomials of a discrete variable on nonuniform lattices\/}
 (in Russian), Keldysh Institute of Applied Mathematics, Preprint
  \#~17, Moscow, 1983.
 
\item{[25]} A.F.~Nikiforov, S.K.~Suslov, V.B.~Uvarov, {\sl Classical
    Orthogonal Polynomials of a Discrete Variable\/}, Springer, Berlin,
    1991.
 
\item {[26]} O.Perron, {\sl Die Lehre von den Kettenbr\"uchen\/},
        $2^{\rm nd}$ edition, Teubner, Leipzig, 1929 = Chelsea,
 
\item{[27]} M.\ Pr\'evost, Dirac masses detection in a density on
    $[-1,1]$ from its moments. Applications to singularities of a
    function, {\sl IMACS Annals on Computing and Appl.\ Math.\/}
   {\bf 9\/} (1991) 365-372.
 
\item{[28]} M.\ Pr\'evost, Dirac masses determination with orthogonal
    polynomials and $\varepsilon-$algorithm. Application to totally
    monotonic sequences, {\sl J.\ Approx.\ Theory\/} {\bf 71\/} (1992)
    175-192.
 
\item{[29]} M.\ Rahman, S.K.\ Suslov, The Pearson equation and the beta
    integral, {\sl SIAM J.\ Math.\ Anal.\ \/} {\bf 25\/} (1994) 646-693.
 
\item{[30]} M.\ Rahman, S.K.\ Suslov, Barnes and Ramanujan-type integrals
    on the $q-$linear lattice,
    {\sl SIAM J.\ Math.\ Anal.\ \/} {\bf 25\/} (1994) 1002-1022.
 
\item{[31]} H.~Stahl, V.~Totik, {\sl General Orthogonal Polynomials\/},
         ({\sl Encyc.\ Math.\ Appl.\/} {\bf 43\/}), Cambridge U.P.,
         Cambridge, 1992.
 
\item{[32]} V.~Totik, Orthogonal polynomials with ratio asymptotics,
    {\sl Proc.\ A.\ M.\ S.} {\bf 114\/} (1992) 491-495.
 
\item{[33]}  W.~Van~Assche, A.P.~Magnus, Sieved orthogonal polynomials and
    discrete measures with jumps dense in an interval, {\sl Proc.\ %
    A.\ M.\ S.\/} {\bf 106\/} (1989) 163-173.

\bye